\input{style/arxiv-general.cfg}
\documentclass[MSNbibl,number,citesort,seceqn,dvips]{arxbj}
\makeatletter
   \@ifpackageloaded{graphicx}{}{\usepackage{graphicx}}
\makeatother

\aid{0}
\volume{21}
\issue{2}
\pubyear{2015}
\firstpage{1231}
\lastpage{1237}
\doi{10.3150/14-BEJ603} 

\makeatletter
\newcommand{\eqref}[1]{(\ref{#1})}
\newtheorem{theorem}{Theorem}[section]
\newtheorem{lemma}[theorem]{Lemma}

\newcommand{\E}{\mathbb{E}}
\renewcommand{\P}{\mathbb{P}}
\newcommand{\R}{\mathbb{R}}
\newcommand{\I}{\mathbb{I}}
\newcommand{\lyg}{\stackrel{\mathrm{def}}{=}}
\newcommand{\vf}{\varphi}
\newcommand{\va}{\mathbf{a}}
\newcommand{\ve}{\varepsilon}
\makeatother

\begin{document}
\begin{frontmatter}

\title{A tight Gaussian bound for weighted sums of Rademacher random variables}
\runtitle{Rademacher--Gaussian}

\begin{aug}
\author[A]{\inits{V.K.}\fnms{Vidmantas Kastytis} \snm{Bentkus}} \and
\author[A]{\inits{D.}\fnms{Dainius} \snm{Dzindzalieta}\corref{}\thanksref{A}\ead[label=e2]{dainiusda@gmail.com}}
\address[A]{Vilnius University Institute of Mathematics and Informatics,
Akademijos 4, Vilnius, Lithuania.\\ \printead{e2}}
\end{aug}

\received{\smonth{2} \syear{2013}}
\revised{\smonth{1} \syear{2014}}

%
\begin{abstract}
Let $\varepsilon_1,\ldots, \varepsilon_n $ be
independent identically distributed Rademacher random variables, that is
$\mathbb{P} \{ \varepsilon_i=\pm1\}=1/2$.
Let
$S_n=a_1\varepsilon_1+\cdots+ a_n\varepsilon_n $,
where $\mathbf{a}=(a_1,\ldots,a_n)\in\mathbb{R}^n$ is a vector such that
${a_1^2+\cdots+a_n^2 \leq1}$.
We find the smallest possible constant $c$ in the inequality
\[
\mathbb{P}\{S_n\geq x\} \leq c \mathbb{P}\{\eta\geq x \} \qquad
\mbox{for all } x\in \mathbb R,
\]
where $\eta\sim N(0,1)$ is a standard normal random variable.
This optimal value
is equal to
\[
c_{\ast}= \bigl(4\mathbb{P}\{\eta\geq\sqrt{2}\} \bigr)^{-1}
\approx3.178.
\]
\end{abstract}

%

\end{frontmatter}

\section{Introduction}

Let $\ve_1,\ldots, \ve_n $ be independent identically distributed
Rademacher random variables, such that
$\P\{ \ve_i=\pm1\}=1/2$. Let
$
S_n=a_1\ve_1+\cdots+ a_n\ve_n $,
where $\va=(a_1,\ldots,a_n)\in\R^n$ is a vector such
that ${a_1^2+\cdots+a_n^2 \leq1}$.

The main result of the paper is the following theorem.

%
\begin{theorem}\label{optx}
Let $\eta\sim N(0,1)$ be a standard
normal random variable. Then, for all $x\in\R$,
%
%
\begin{equation}
\P\{S_n\geq x\} \leq c \P\{\eta\geq x \},
\end{equation}
with the constant $c$ equal to
\[
c_{\ast}:= \bigl(4\P\{\eta\geq\sqrt{2}\} \bigr)^{-1}
\approx3.178.
\]
\end{theorem}

The value $c=c_\ast$ is the best possible since
\eqref{optx}
becomes equality if $n\geq2$, $x=\sqrt{2}$ and $S_n=(\ve_1+\ve
_2)/\sqrt{2}$.

Inequality \eqref{optx} was first obtained by Pinelis \cite
{pinelis94} with $c\approx4.46$. Bobkov, G\" otze and Houdr\'e (BGH)
\cite{bobkov2001gaussian} gave a simple proof of \eqref{optx} with
constant factor ${c\approx12.01}$. Their method was to use induction
on $n$ together with the inequality
%
%
\begin{equation}
\label{kkk} \tfrac{1}{2}\P\{\eta\geq A\}+\tfrac{1}{2}\P\{\eta\geq B
\} \leq\P \{\eta \geq x\} \qquad\mbox{for all } x\geq\sqrt{3} \mbox{ and } \tau
\in[0,1],
\end{equation}
where $A:=\frac{x-\tau}{\sqrt{1-\tau^2}}$ and $B:=\frac{x+\tau
}{\sqrt{1-\tau^2}}$. Using a method similar to the one in BGH \cite
{bobkov2001gaussian}, Bentkus \cite{bentkus2007measure} proved \eqref
{optx} with $c\approx4.00$ and conjectured that the optimal constant
in \eqref{optx} is $c_\ast$. Further progress was achieved by Pinelis
\cite{pinelis2007}, where \eqref{optx} was proved with $c\approx1.01
c_{\ast}$.

Let us briefly outline our strategy of the proof. For $x\leq\sqrt{2}$
Theorem~\ref{optx} follows from the symmetry of $S_n$. For $x\geq
\sqrt {2}$ we consider two cases separately. If $x\in (\sqrt{2},\sqrt {3} )$ and all $a_i$'s are ``small'' we use Berry--Esseen
inequality. Otherwise we use induction on $n$ together with Chebyshev
type inequality presented in Lemma~\ref{inequality}. We remark that the
analysis of weighted sums of random variables based on separate study
of these two cases has proved recently to be effective idea, see \cite
{rudelson2009smallest}. In \cite{pinelis2012}, this idea was used to
obtain asymptotically Gaussian bound
\[
\P\{S_n\geq x\}\leq\P\{\eta\geq x \} \biggl(1+\frac{C}{x}
\biggr),
\]
where $C\approx14.10\ldots$\,.

A standard application of bounds like
\eqref{optx}, following Efron \cite{efron69}, is to the Student's
statistic and to
self-normalized sums. For example,
if random variables ${ {X}_{1},\ldots, {X}_{n} }$
are independent (not necessary identically distributed), symmetric and
not all identically equal to zero, then
the statistic
\[
T_n=({{X}_{1}+\cdots+ {X}_{n}})/\sqrt
{X_1^2+\cdots+X_n^2}
\]
is
sub-Gaussian and
%
%
\begin{equation}
\label{optselfx} \P\{ T_n \geq x\} \leq c_\ast \P\{\eta\geq x
\} \qquad\mbox{for all }x\in\R.
\end{equation}
The latter inequality is optimal since it turns into an equality
if $n=2$, $x=\sqrt{2}$ and $X_1=\ve_1$, $X_2=\ve_2$. This inequality
was previously obtained in \cite{pinelis94,pinelis2007} with
constants $4.46$ and ${\approx}1.01c_{\ast}$ in place of $c_{\ast}$.

\section{Proofs}
\label{sectionm}

In this section, we use the following notation
%
%
\begin{equation}
\label{cheb2} \tau=a_1, \qquad\vartheta= \sqrt{1-
\tau^2}, \qquad I(x)=\P\{\eta\geq x\},\qquad\varphi(x)=-I'(x),
\end{equation}
that is, $I(x)$ is the tail probability for standard normal random
variable $\eta$ and $\varphi(x)$ is the standard normal density.
Without loss of generality, we assume that
$ {a_1^2+\cdots+a_n^2 = 1}$ and $ {a_1\geq\cdots\geq a_n\geq0}$.
Using \eqref{cheb2} we have $S_n =\tau\ve_1+ \vartheta X$ with $ X=
(a_2\ve_2+\cdots+a_n\ve_n )/ \vartheta$.
The random variable $X$ is symmetric and independent of $\ve_1$.
It is easy to check that
$\E X^2=1$ and
%
%
\begin{equation}
\label{cheb3} \P\{S_n \geq x\} = \tfrac{1}{2} \P\{X \geq A\}+
\tfrac{1}{2} \P\{X \geq B\},
\end{equation}
where $A=\frac{x-\tau}\vartheta$ and $
B= \frac{x+\tau}\vartheta$.

We start with a simple Chebyshev type inequality.

%
\begin{lemma}\label{inequality}
Let $s>0$ and $0\leq a \leq b$. Then, for any random variable $Y$ we have
%
%
\begin{equation}
\label{in1} a^s\P \bigl \{|Y|\geq a \bigr \}+ \bigl(b^s-a^s
\bigr)\P\bigl \{|Y|\geq b\bigr \}\leq\E|Y|^s.
\end{equation}
If $Y$ is symmetric, then
%
%
\begin{equation}
\label{in2} a^s\P\{Y\geq a\}+ \bigl(b^s-a^s
\bigr)\P\{Y\geq b\}\leq\E|Y|^s/2.
\end{equation}
\end{lemma}

\begin{pf}
It is clear that
\eqref{in1}
implies
\eqref{in2}. To prove \eqref{in1}, we use the obvious inequality
%
%
\begin{equation}
\label{in3} a^s\I\bigl \{|Y|\geq a\bigr \}+ \bigl(b^s-a^s
\bigr)\I\bigl \{|Y| \geq b\bigr \}\leq|Y|^s,
\end{equation}
where $\I\{E\}$ stands for the indicator function of the event $E$.
Taking expectation, we get
\eqref{in1}. 
\end{pf}

In similarity to \eqref{in1}, one can derive a number
of inequalities stronger than the standard Chebyshev inequality
$\P\{S_n\geq x\}\leq1/(2x^2)$. For example, instead of
$ {\P\{S_n\geq1\}\leq1/2}$ we have the much stronger
\[
\P\{S_n\geq1\}+\P\{S_n\geq\sqrt{2}\}+\P\{S_n
\geq\sqrt{3}\} +\cdots \leq1/2.
\]

We will make use of Lyapunov type bounds with explicit constants
for the remainder term in the Central limit theorem.
Let ${ {X}_1, {X}_2,\ldots} $ be independent random variables such that
$\E X_j=0$ for all $j$.
Set $\beta_j = \E|X_j|^3$.
Assume that the sum ${Z={ {X}_1+ {X}_2+\cdots}}$ has unit variance.
Then there exists an absolute constant,
say $c_{\mathrm{L}}$, such that
%
%
\begin{equation}
\label{clt0a} \bigl |\P\{ Z\geq x\}- I(x)\bigr |\leq c_{\mathrm{L}} (
\beta_1+ \beta_2+\cdots).
\end{equation}
It is known that
$c_{\mathrm{L}} \leq0.56\ldots$ \cite
{tyurin2012refinement,shevts2010}. Note that we actually do not need
the best known bound for $c_\mathrm{L}$. Even $c_\mathrm{L}=0.958$ suffices to prove
Theorem~\ref{optx}.

Replacing $X_j$ by $a_j\ve_j$ and using $\beta_j\leq\tau a^2_j$ for
all $j$, the
inequality \eqref{clt0a} implies
%
%
\begin{equation}
\label{clt} \bigl |\P\{ S_n\geq x\}- I(x)\bigr |\leq c_{\mathrm{L}}\tau.
\end{equation}

\begin{pf*}{Proof of Theorem~\ref{optx}}
For $x\leq\sqrt{2}$ Theorem~\ref{optx} follows from the symmetry of $S_n$
and Chebyshev's inequality (first it was implicitly shown in \cite
{bentkus2007measure}, later in \cite{pinelis2007}).
In the case $x\geq\sqrt{2}$, we argue by induction on $n$.
However, let us first provide a proof of Theorem~\ref{optx}
in some special cases where induction fails.

Using the bound \eqref{clt}, let us prove Theorem~\ref{optx} under the
assumption that
%
%
\begin{equation}
\label{clt0} \tau\leq\tau_{\mathrm{L}} \lyg (c_{\ast}-1) I(\sqrt
{3})/c_{\mathrm{L}} \quad\mbox{and}\quad x\leq\sqrt{3}.
\end{equation}
Using $c_{\mathrm{L}} = 0.56$, the numerical value of
$\tau_{\mathrm{L}}$ is
$0.16\ldots$\,. In order to prove Theorem~\ref{optx}
under the assumption \eqref{clt0}, note that the
inequality \eqref{clt} yields
%
%
\begin{equation}
\label{clt1} \P\{ S_n\geq x\}\leq I(x)+ \tau c_{\mathrm{L}}.
\end{equation}
If\vspace*{2pt} the inequality \eqref{clt0} holds,
the right-hand side of \eqref{clt1}
is clearly bounded from above by ${ c_\ast
I(x)}$ for $x\leq\sqrt{3}$.

For $x$ and $\tau$ such that \eqref{clt0} does not hold we use
induction on $n$.
If $n=1$, then we have $S_n=\ve_1$ and
Theorem~\ref{optx} is equivalent to the trivial inequality $1/2 \leq
c_\ast I(1)$.

Let us assume that Theorem~\ref{optx} holds for $n\leq k-1$ and prove it
for $n=k$.

Firstly we consider the case $x\geq\sqrt{3}$.
We replace $ S_n$ by $S_k$
with
$ X= (a_2\ve_2+\cdots+a_k\ve_k)/ \vartheta$ in \eqref{cheb3}.
We can estimate
the latter two probabilities in \eqref{cheb3}
applying the induction hypothesis
$\P\{X \geq y\}\leq c_\ast I(y)$. We get
%
%
\begin{equation}
\label{cheb4} 
\P\{S_k \geq x\} \leq c_\ast
I(A)/2+c_\ast I(B)/2.
\end{equation}
In order to conclude the proof, it suffices
to show that the right-hand side of \eqref{cheb4}
is bounded from above by $c_\ast I(x)$, that is,
that the inequality $ I(A)+ I(B)\leq2I(x)$ holds.
As $x\geq\sqrt{3}$ this follows by the inequality \eqref{kkk}.

In the remaining part of the proof, we can assume that $x\in(\sqrt {2},\sqrt{3})$ and $\tau\geq\tau_{\mathrm{L}}$.
In this case in order to prove Theorem~\ref{optx}, we have to improve the
arguments used to estimate
the right-hand side of \eqref{cheb3}.
This is achieved by applying the Chebyshev type inequalities of
Lemma~\ref{inequality}. By Lemma~\ref{inequality},
for any symmetric $X$ such that
$\E X^2=1$, and $0\leq A\leq B$, we have
%
%
\begin{equation}
\label{in4a} A^2\P\{X\geq A\}+ \bigl(B^2-A^2
\bigr) \P\{X\geq B\} \leq1/2.
\end{equation}
By \eqref{cheb2},
we can rewrite \eqref{in4a} as
%
%
\begin{equation}
\label{in4a1} (x-\tau)^2\P\{X\geq A\}+4x\tau\P\{X\geq B\} \leq
\vartheta^2/2.
\end{equation}

For $x\in(\sqrt{2},\sqrt{3})$ and $\tau\geq\tau_{\mathrm{L}}$ we
consider the cases
\[
\mbox{(i)}\ (x-\tau)^2\geq4x\tau\quad\mbox{and}\quad \mbox{(ii)}\
(x-\tau)^2\leq4x\tau
\]
separately. We denote the sets of points $(x,\tau)$ such that $x\in
(\sqrt{2},\sqrt{3})$, $\tau\geq\tau_{\mathrm{L}}$ and (i) or (ii)
holds by $E_1$ and $E_2$, respectively.

(i) Using \eqref{cheb3}, \eqref{in4a1} and the induction hypothesis
we get
%
%
\begin{equation}
\label{in6d} \P\{S_k \geq x\}\leq\frac{D
\P\{X\geq B\}+\vartheta^2/2}{2(x-\tau)^2} \leq
\frac{c_\ast D
I( B)+\vartheta^2/2}{2(x-\tau)^2},
\end{equation}
where $X=(a_2\ve_2+\cdots+a_k\ve_k)/\vartheta$ and $D=(x-\tau
)^2-4x\tau$.

In order to finish the proof of Theorem~\ref{optx}
(in this case) it suffices
to show that the right-hand side of \eqref{in6d}
is bounded above by $c_\ast I( x)$. In other words,
we have to check that the function
%
%
\begin{equation}
\label{in7e} f(x,\tau)\equiv f \lyg { \bigl((x-\tau)^2-4x\tau
\bigr) c_\ast I(B)-2c_\ast(x-\tau)^2I(x)+
\vartheta^2/2},
\end{equation}
is negative on $E_1$, where $B=(x+\tau)/\vartheta$.

By Lemma~\ref{lem1} below, we have
%
%
\begin{equation}
\label{ineq3} f(x,\tau)\leq f(\sqrt{3},\tau)=:g(\tau).
\end{equation}
Since $\tau\leq(3-2\sqrt{2})x$
the inequality ${f\leq0}$ on $E_1$ follows from Lemma~\ref{lemg}, below.

(ii) Using \eqref{cheb3}, \eqref{in4a1} and induction hypothesis we get
%
%
\begin{equation}
\label{in6db} \P\{S_k \geq x\} \leq\frac
{{C
\P\{X\geq A\}+\vartheta^2/2}} {
8x\tau}\leq
\frac{ C/(2A^2)+\vartheta^2/2}{8x\tau},
\end{equation}
where $X=(a_2\ve_2+\cdots+a_k\ve_k)/\vartheta$ and $C=4x\tau
-(x-\tau)^2$.

In order to finish the proof (in this case) it suffices
to show that the right-hand side of \eqref{in6db}
is bounded above by $c_\ast I( x)$. In other words,
we have to check that
%
%
\begin{equation}
\label{in6db1} C/ \bigl(2A^2 \bigr)+\vartheta^2/2
\leq{8x \tau}c_\ast I(x) \qquad\mbox{on } E_2.
\end{equation}
Recalling that
$C=4x\tau-(x-\tau)^2$, $A =(x-\tau)/\vartheta$,
inequality \eqref{in6db1}
is equivalent to
%
%
\begin{equation}
\label{in11} h \lyg \frac{1-\tau^2}{(x-\tau)^2}- 4c_\ast I(x)\leq0 \qquad
\mbox{on } E_2.
\end{equation}
Inequality \eqref{in11} follows from Lemma~\ref{lemh}, below.
The proof of Theorem~\ref{optx} is complete.
\end{pf*}
%

%
\begin{lemma}\label{lemg}
The function $g$ defined by \eqref{ineq3} is negative for all $\tau
\in[\tau_{\mathrm{L}},(3-2\sqrt{2})\sqrt{3}]$.
\end{lemma}

%
%
\begin{lemma}\label{lem2}
$I'(B)\geq\vartheta I'(x)$ on $E_1$.
\end{lemma}

%
%
\begin{lemma}\label{lem3}
$I(B)\geq I(x)+I'(x)\tau$ on $E_1$.
\end{lemma}

%
%
\begin{lemma}\label{lem1}
The partial derivative $\partial_x f$ of the function $f$ defined by
\eqref{in7e} is positive on $E_1$.
\end{lemma}

%
%
\begin{lemma}\label{lemh}
The function $h$ defined by \eqref{in11} is negative on $E_2$.
\end{lemma}

\begin{pf*}{Proof of Lemma~\ref{lemg}}
Since $g(\tau_{\mathrm{L}})<0$ it is sufficient to show that $g$ is a
decreasing
function for $\tau_{\mathrm{L}}\leq\tau\leq(3-2\sqrt{2})\sqrt{3}$.
Note that
\[
g(\tau) = { \bigl((\sqrt{3}-\tau)^2-4\sqrt{3}\tau \bigr)
c_\ast I(B)+ \bigl(1-\tau^2 \bigr)/2}- {2c_\ast (
\sqrt{3}-\tau)^2}I(\sqrt{3})
\]
and
\begin{eqnarray*}
g'(\tau)&=& (2\tau-6\sqrt{3}) c_\ast I(B) - \bigl((
\sqrt{3}-\tau)^2-4\sqrt{3}\tau \bigr) c_\ast \vf(B) (1+\tau
\sqrt{3})\vartheta^{-3}
\\
&&{}-\tau+ {4c_\ast (\sqrt{3}-\tau)}I(
\sqrt{3}),
\end{eqnarray*}
where $\vf$ is the standard normal distribution.
Hence
\[
g'(\tau)\leq w(\tau) \lyg (2\tau-6\sqrt{3}) c_\ast I(B) -
\tau+ {4c_\ast (\sqrt{3}-\tau)}I(\sqrt{3}).
\]
Note that the value of $B$ in previous three displayed formulas should
also be computed with $x=\sqrt{3}$.
Using Lemma~\ref{lem3}, we get
\[
g'(\tau)\leq -2c_\ast (\sqrt{3}+\tau) I(\sqrt{3}) +
2c_\ast \tau(3\sqrt{3}-\tau) \vf(\sqrt{3})-\tau\lyg Q(\tau)
\]
with
\[
Q(\tau)=-\alpha\tau^2+\beta\tau- \gamma,\qquad \alpha=0.56\ldots,
\qquad\beta=1.67\ldots, \qquad\gamma =0.45\ldots.
\]

Clearly, $Q$ is negative on the interval $[\tau_{\mathrm
{L}},(3-2\sqrt {2})\sqrt{3} ]$. It follows that $g'$ is negative, and $g$
is decreasing on $[\tau_{\mathrm{L}},(3-2\sqrt{2})\sqrt{3}]$.
\end{pf*}

\begin{pf*}{Proof of Lemma~\ref{lem2}}
Since ${I'=-\vf}$ by \eqref{cheb2}, the inequality $ { I'(B)\geq
\vartheta I'(x)}$ is equivalent to
%
%
\[
u(\tau) \lyg { \bigl(1-\tau^2 \bigr)} \exp \biggl\{
\frac{(x+\tau)^2}{1-\tau^2}-x^2 \biggr\} -1\geq0.
\]
Since $u(0)=0$, it suffices to check that
$u'\geq0$. Elementary calculations show that ${u'\geq0}$
is equivalent to the trivial inequality
$x+\tau^2 x +\tau x^2 +\tau^3\geq0$. 
\end{pf*}

\begin{pf*}{Proof of Lemma~\ref{lem3}}
Let $g(\tau)=I(B)$. Then the
inequality ${I(B)\geq I(x)+I'(x)\tau}$
turns into $g(\tau)\geq g(0)+g'(0)\tau$. The latter inequality holds
provided that ${g''(\tau)\geq0}$.
Next, it is easy to see that
$g'(\tau)=-\vf(B)B'$ and $ {g''(\tau) = (BB^{\prime2}-B'')\vf(B)}$.
Hence, to verify that $g''(\tau)\geq0$ we verify that
$ BB^{\prime2}-B''\geq0$.
This last inequality is equivalent to
$ -2 +2x^2 +x^3\tau+x^2\tau^2+x\tau+2x\tau^3+3\tau^2\geq0$,
which holds since $x\geq1$. The proof of Lemma~\ref{lem3}
is complete.
\end{pf*}
%

\begin{pf*}{Proof of Lemma~\ref{lem1}}
We have
\[
\partial_x f=2 (x-3\tau ) c_\ast I(B)+ D c_\ast
I'(B)/\vartheta-4c_\ast(x-\tau)I(x)- 2c_\ast(x-
\tau)^2I'(x).
\]
We have to show that $\partial_x f\geq0$ on $E_1$.
Using Lemma~\ref{lem2},
we can reduce this to the inequality
%
%
\begin{equation}
\label{der1} 2 (x-3\tau ) I(B)-(x+\tau)^2 I'(x)-4(x-
\tau)I(x)\geq0.
\end{equation}
On $E_1$ we have that $0\leq\tau\leq(3-2\sqrt{2})x$, so $x-3\tau
\geq x-3(3-2\sqrt{2})x=(6\sqrt{2}-8)x>0$.
By Lemma~\ref{lem3} we have that left-hand side of \eqref{der1} is
bigger than
\begin{eqnarray*}
&& 2(x-3\tau) \bigl(I(x)+I'(x)\tau \bigr)-(x+\tau)I'(x)-4(x-
\tau)I(x)
\\
&&\quad= -2(x+\tau)I(x)- \bigl(x^2+7\tau^2
\bigr)I'(x).
\end{eqnarray*}
Inequality \eqref{der1} follows by the inequality
$-(x^2+7\tau^2)I'(x)\geq\alpha x(x+\tau)\vf(x)>2(x+\tau)I(x)$ on
$E_1$ with $\alpha=4\sqrt{14}-14$, where the second inequality
follows from the fact that $ {\vf(x)x/I(x)}$ increases for $x>0$ and
is larger than $2/\alpha$ for $x=\sqrt{2}$. The proof of Lemma~\ref
{lem1} is complete.
\end{pf*}

\begin{pf*}{Proof of Lemma~\ref{lemh}}
It is easy to see that the function $h$
attains its maximal value at $\tau=1/x$. Hence,
it suffices to check \eqref{in11} with $\tau=1/x$, that is, that
for $\sqrt{2}\leq x\leq\sqrt{3}$ the inequality
$g(x) \lyg 1-
4c_\ast(x^2-1) I(x)\leq0
$ holds.
Using\vspace*{2pt} $4 c_\ast I(\sqrt{2})=1$, we
have $g(\sqrt{2})=0$ and $g(\sqrt{3})<0$. Next,
$g'(x)=-8c_\ast xI(x)+4c_\ast(x^2-1)\vf(x)$, so $g'(\sqrt{2})<0$ and
$g'(\sqrt{3})>0$.
We have\vspace*{2pt} that $g''(x)=4c_{\ast} ((5-x^2)x\vf(x)-2I(x) )$.
Since $I(x)\leq\vf(x)/x$ we have that $g''(x)\geq4c_{\ast}
((5-x^2)x\vf(x)-2\vf(x)/x )=4c_{\ast}\vf
(x)/x((5-x^2)x^2-2)\geq8c_\ast\vf(x)/x>0$\vspace*{2pt} for $x\in(\sqrt{2},\sqrt {3})$. The proof of Lemma~\ref{lemh} is complete.
\end{pf*}
%

\section*{Acknowledgements}

This research was funded by a grant (No. MIP-47/2010) from the Research
Council of Lithuania.


%

\printhistory


\begin{thebibliography}{9}

\bibitem{bentkus2007measure}
%
\begin{barticle}[mr]
\bauthor{\bsnm{Bentkus},~\bfnm{Vidmantas}\binits{V.}}
(\byear{2007}).
\btitle{On measure concentration for separately {L}ipschitz functions
in product spaces}.
\bjournal{Israel J. Math.}
\bvolume{158}
\bpages{1--17}.
\bid{doi={10.1007/s11856-007-0001-2}, issn={0021-2172}, mr={2342455}}
\end{barticle}
%
\bptok{imsref}%
\endbibitem

\bibitem{bobkov2001gaussian}
%
\begin{barticle}[mr]
\bauthor{\bsnm{Bobkov},~\bfnm{Sergey~G.}\binits{S.G.}},
\bauthor{\bsnm{G{\"o}tze},~\bfnm{Friedrich}\binits{F.}} \AND
\bauthor{\bsnm{Houdr{\'e}},~\bfnm{Christian}\binits{C.}}
(\byear{2001}).
\btitle{On {G}aussian and {B}ernoulli covariance representations}.
\bjournal{Bernoulli}
\bvolume{7}
\bpages{439--451}.
\bid{doi={10.2307/3318495}, issn={1350-7265}, mr={1836739}}
\end{barticle}
%
\bptok{imsref}%
\endbibitem

\bibitem{efron69}
%
\begin{barticle}[mr]
\bauthor{\bsnm{Efron},~\bfnm{Bradley}\binits{B.}}
(\byear{1969}).
\btitle{Student's {$t$}-test under symmetry conditions}.
\bjournal{J. Amer. Statist. Assoc.}
\bvolume{64}
\bpages{1278--1302}.
\bid{issn={0162-1459}, mr={0251826}}
\end{barticle}
%
\bptok{imsref}%
\endbibitem

\bibitem{pinelis94}
%
\begin{barticle}[mr]
\bauthor{\bsnm{Pinelis},~\bfnm{Iosif}\binits{I.}}
(\byear{1994}).
\btitle{Extremal probabilistic problems and {H}otelling's {$T\sp2$}
test under a symmetry condition}.
\bjournal{Ann. Statist.}
\bvolume{22}
\bpages{357--368}.
\bid{doi={10.1214/aos/1176325373}, issn={0090-5364}, mr={1272088}}
\end{barticle}
%
\bptok{imsref}%
\endbibitem

\bibitem{pinelis2007}
%
\begin{barticle}[mr]
\bauthor{\bsnm{Pinelis},~\bfnm{Iosif}\binits{I.}}
(\byear{2007}).
\btitle{Toward the best constant factor for the
{R}ademacher--{G}aussian tail comparison}.
\bjournal{ESAIM Probab. Stat.}
\bvolume{11}
\bpages{412--426}.
\bid{doi={10.1051/ps:2007027}, issn={1292-8100}, mr={2339301}}
\end{barticle}
%
\bptok{imsref}%
\endbibitem

\bibitem{pinelis2012}
%
\begin{barticle}[mr]
\bauthor{\bsnm{Pinelis},~\bfnm{Iosif}\binits{I.}}
(\byear{2012}).
\btitle{An asymptotically {G}aussian bound on the {R}ademacher tails}.
\bjournal{Electron. J. Probab.}
\bvolume{17}
\bpages{no. 35, 22}.
\bid{doi={10.1214/EJP.v17-2026}, issn={1083-6489}, mr={2924368}}
\end{barticle}
%
\bptok{imsref}%
\endbibitem

\bibitem{rudelson2009smallest}
%
\begin{barticle}[mr]
\bauthor{\bsnm{Rudelson},~\bfnm{Mark}\binits{M.}} \AND
\bauthor{\bsnm{Vershynin},~\bfnm{Roman}\binits{R.}}
(\byear{2009}).
\btitle{Smallest singular value of a random rectangular matrix}.
\bjournal{Comm. Pure Appl. Math.}
\bvolume{62}
\bpages{1707--1739}.
\bid{doi={10.1002/cpa.20294}, issn={0010-3640}, mr={2569075}}
\end{barticle}
%
\bptok{imsref}%
\endbibitem

\bibitem{shevts2010}
%
\begin{barticle}[mr]
\bauthor{\bsnm{Shevtsova},~\bfnm{I.~G.}\binits{I.G.}}
(\byear{2010}).
\btitle{Refinement of estimates for the rate of convergence in
{L}yapunov's theorem}.
\bjournal{Dokl. Akad. Nauk}
\bvolume{435}
\bpages{26--28}.
\bid{doi={10.1134/S1064562410060062}, issn={0869-5652}, mr={2790498}}
\end{barticle}
%
\bptok{imsref}%
\endbibitem

\bibitem{tyurin2012refinement}
%
\begin{barticle}[mr]
\bauthor{\bsnm{Tyurin},~\bfnm{I.~S.}\binits{I.S.}}
(\byear{2011}).
\btitle{Improvement of the remainder in the {L}yapunov theorem}.
\bjournal{Teor. Veroyatn. Primen.}
\bvolume{56}
\bpages{808--811}.
\bid{issn={0040-361X}, mr={3137072}}
\bptnote{check year}%
\end{barticle}
%
\bptok{imsref}%
\endbibitem

\end{thebibliography}
\end{document}